\documentclass[11pt]{amsart}
\usepackage{amsmath}
\usepackage{palatino}
  \usepackage{paralist}
  \usepackage[misc]{ifsym}
  \usepackage{graphics} 
  \usepackage{epsfig} 

\usepackage{color}

\newtheorem{theorem}{Theorem}[section]
\newtheorem{corollary}{Corollary}

\newtheorem{lemma}[theorem]{Lemma}
\newtheorem{proposition}{Proposition}

\theoremstyle{definition}

\newtheorem{remark}{Remark}

\newcommand{\R}{\mathbb{R}}
\newcommand{\C}{\mathbb{C}}

\newcommand{\be}{\begin{equation}}
\newcommand{\ee}{\end{equation}}

\title[ Klein-Gordon Schr\"{o}dinger system]
      {Global Attractor for two-dimensional weakly damped nonlinear parametric Klein-Gordon Schr\"{o}dinger system}

\author[Olivier Goubet and Marilena Poulou]{}


 \email{olivier.goubet@univ-lille.fr}
 \email{  mpoulou@uniwa.gr}

\begin{document}
\maketitle

\centerline{\scshape Olivier Goubet}
\medskip
{\footnotesize
\centerline{Laboratoire Paul Painlevé CNRS UMR 8524,}

\centerline{Université de Lille, INRIA,}
\centerline{59 655 Villeneuve d'Ascq cedex.}

\medskip
\centerline{and }

} 

\medskip

\centerline{\scshape  Marilena Poulou}
{\footnotesize
\centerline{Department of Mechanical Engineering,}

\centerline{University of West Attica,}
\centerline{Egaleo, 12241, Athens.}

\bigskip


\begin{abstract}
We consider a damped forced nonlinear Klein-Gordon Schr\"{o}dinger (KGS) system
known as the parametric KGS system with Yukawa coupling where the space variable belongs to $\R^2$.
We prove the existence of a global attractor in dimension $2$.
Moreover we prove that this global attractor is regular and has finite fractal dimension.
\end{abstract}

\section{Introduction}

We consider a Klein-Gordon Schr\"{o}dinger system (KGS) with Yukawa coupling.
This system describes the nonlinear interaction between high frequency
electron waves and low frequency ion plasma waves in a homogeneous magnetic field,
adapted to model the upper hybrid plasma
heating scheme (see \cite{KSX} and the references therein).
This system couples a Schr\"odinger equation with
a wave equations and reads
\begin{eqnarray}
 \label{e1.1} \psi_t +i\Delta \psi  +\alpha\psi &= & i\phi \psi +f, \\
 \label{e1.2}
 \phi_{tt}-\Delta \phi +\phi+\gamma\phi_t  &= &  -Re \psi_{x_1} +g.
\end{eqnarray}

\noindent Here the unknowns $\psi$ and $\phi$ map respectively
$\R_t\times \R^2_x$ into respectively $\C$ and $\R$. The positive parameters $\alpha$
and $\gamma$ are the damping parameters that are supposed to be small, i.e. $\alpha, \gamma<<1$. In the right hand side $f, g$ are external forces
which do not depend on $t$ and that belong to $L^2(\R^2)$.

This weakly damped forced KGS system partakes of infinite-dimensional
dynamical systems  \cite{Hale}, \cite{lady}, \cite{MZ},\cite{Raugel}, \cite{Temam}.
This system was studied in $H^2\times H^2\times H^1$ in the case
of a bounded domain and homogneous Dirichlet conditions in \cite{PS1} and in the unbounded domain in \cite{PS2}. In both works the authors proved the existence of a global attractor that attracts all bounded sets in the norm topology.

We plan to address the issue of the existence of a global attractor
in the energy space $\mathbb{E}_1=H^1(\R^2)\times H^1(\R^2) \times L^2(\R^2)$.
This issue was studied in dimension $1$ in \cite{GoubetPoulou2025}.

Our first result states as follows

\begin{theorem}\label{th1}
There exists a compact global attractor $\mathcal{A}$ in $\mathbb{E}_1$.
This compact global attractor is moreover a compact subset of $H^2(\R^2)\times H^2(\R^2) \times H^2(\R^2)$.
\end{theorem}

This result is sharp since the global attractor contains the stationary solutions of the equation, i.e.
solution $(\psi^*,\phi^*,0)$ that solves
\begin{eqnarray}
 \label{e.1} i\Delta \psi^*  +\alpha\psi^* &= & i\phi^* \psi^* +f, \\
 \label{e.2}
 -\Delta \phi^* +\phi^* &= &  -Re \psi_{x_1}^* +g.
\end{eqnarray}

\noindent These stationary solutions belong to $H^2(\R^2 )\times H^2(\R^2 )\times \{0\}.$

Our second result asserts

\begin{theorem}\label{th2}
The fractal dimension of $\mathcal{A}$ is finite.
\end{theorem}

\begin{remark}
We point out that we do not need any decay assumption on $f,g$.
\end{remark}

The remaining of the article is organized as follows. In a Section 2
we describe the mathematical framework and recall some functional analysis
results. In Section 3 we describe some results concerning the initial value problem.
In Section 4 and 5 we prove the existence of a global attractor in the energy space.
In section 6 we prove the regularity of the global attractor.
In section 7 we prove that this attractor has finite dimension.
\smallskip

To complete this introduction, let us introduce some notations. We set $c$ for a numerical constant that may change from one line to one
another. The letter $K$ denotes some generic constant that
depends on the data $\alpha, \gamma, f, g$.

\section{Mathematical framework}

\subsection{Some inequalities}

The Hilbert space $L^2(\R^2)$ is equipped with the scalar product
 $$\langle u,v\rangle= {\rm Re}\int_{\mathbb{R}^2}u(x) \overline{v}(x)dx.$$
 We also use standard Sobolev and Banach spaces as $H^1(\R^2)$ or $L^\infty(\R^2)$.

In dimension $2$, $H^1(\R^2)$ is not a Banach algebra. Besides, for any function
in $H^2(\R^2)$ we have the so-called Agmon inequality

\begin{equation}\label{agmon}
||u||^2_{L^\infty(\R^2)}\leq c ||u||_{L^2(\R^2)}||\Delta u||_{L^2(\R^2)}.
\end{equation}

\noindent Another useful inequality is the Gagliardo-Nirenberg inequality
that reads for any function in $H^1(\R^2)$

\begin{equation}\label{gn}
||u||^2_{L^4(\R^2)}\leq c ||u||_{L^2(\R^2)}||\nabla u||_{L^2(\R^2)}.
\end{equation}

We now introduce the so-called Strichartz estimates.
Let us recall the following statement (see \cite{Cazenave}).
Set $\psi(t)=S(t)\psi_0$ for the solution of the free
Schr\"odinger equation
$$ \psi_t+i\Delta \psi=0, $$
\noindent with $\psi(0)=\psi_0$. Then

\begin{lemma}\label{strichartz}
We define an admissible pair $(q,r)$ as $\frac 2 q + \frac 2 r=1$ and $(q,r)\neq (2,\infty)$.
Let $(\tilde q',\tilde r')$ to be the conjugate exponents of another admissible pair. Then
we have the following space-time estimates
$$ ||S(t)\psi_0||_{L^q_tL^r_x}\leq c ||\psi_0||_{L^2(\R^2)},$$
$$ ||\int_0^t S(t-s) F(s)ds||_{L^q_tL^r_x}\leq c ||F||_{L^{\tilde q'}_tL^{\tilde r'}_x}. $$
\end{lemma}

\subsection{A memory operator}

\noindent Consider the equation \eqref{e1.2} that is linear.
 Following \cite{Bourgain}, we introduce the new variable $v=\varphi+i(Id-\Delta)^{-\frac12}\varphi_t$
that is solution to the Schr\"odinger like equation

\begin{equation}\label{nls}
v_t+i\gamma {\rm Im} v+i(Id-\Delta)^\frac12 v=i (Id-\Delta )^{-\frac12}(-Re \psi_{x_1} +g) .
\end{equation}

We use the following statement (see \cite{GoubetPoulou2025} for a proof)

 \begin{lemma}\label{lemma} Let $\sigma$ in $\R$ and
a function $h(t,x)$ that is continuous and bounded from $\R_t$
into $H^{\sigma}(\R^2)$. Consider $(u,u_t)$ the solution to
\be\label{wave}
u_{tt}+\gamma u_t-\Delta u+u=h,
\ee
\noindent that is continuous and bounded from $\R_t$ into $H^\sigma(\R^2)\times H^{\sigma-1}(\R^2)$.
Then $(u,u_t)$ is a continuous and bounded function from $\R_t$ into $H^{\sigma+1}(\R^2)\times H^{\sigma}(\R^2)$.
\end{lemma}

This a mere consequence of the following representation formula.
Consider the family of linear operators defined by $U_\gamma(t)v_0=v(t)$ where $v$ solves
\eqref{nls} with $h=0$, and with $v(0)=v_0$. Then

\begin{equation}\label{rep}
v(t)=i\int_{-\infty}^t U_{\gamma}(t-s) (Id-\Delta)^{-\frac12} h(s)ds.
\end{equation}

\noindent The smoothing comes from the operator $(1-\Delta)^{-\frac12}$
and the convergence of the integral in $-\infty$ comes from the damping term $i\gamma {\rm Im} v$.

\subsection{Some compactness and regularity lemma}

We first recall without proof the following lemma

\begin{lemma}\label{compact}
Consider a closed set $\mathcal{E}$ in Banach space $E$ such that for any $\varepsilon>0$
there exists a compact set $K_\varepsilon$ of $E$ such that $\mathcal{E}\subset K_\varepsilon+B_E(0,\varepsilon)$,
where $B_E(0,\varepsilon)$ is the closed ball in $E$ centered at $0$ and of radius
$\varepsilon$. Then $\mathcal{E}$ is a compact subset of $E$.
\end{lemma}

This lemma is useful to overcome the difficulty that the embedding $H^1(\R^2)\subset L^2(\R^2)$
is not compact. If we have a bounded subset $\mathcal{B}$ of $H^1(\R^2)$ that satisfies a uniform tail estimate
(see Lemma \ref{luwang} below), then since by Rellich's theorem $H^1(\R^2)\subset L^2([-R,R]^2)$
is compact, then $\mathcal{B}$ is a compact subset of $L^2(\R^2)$.

We complete this section with the statement

\begin{lemma}\label{tgtbt}
Consider a function $w$ in $H^1(\R^2)$ such that $w=u^m+v^m$ where $u^m$
is a bounded sequence in $H^2(\R^2)$ and $v^m$ converges to $0$ when
$m$ diverges to $+\infty$, then $w$ belongs to $H^2(\R^2)$.
\end{lemma}

\begin{proof}
Introduce the projection operator, for $\mathcal{F}$ the Fourier transform
and $\chi(\xi)=1$ if $\xi \in [-N,N]^2$ and $\chi(\xi)=0$ if not

\begin{equation}\label{pn}
P_Nw= \mathcal{F}^{-1}\left( \chi \mathcal{F}(u) \right).
\end{equation}

\noindent Then by Bernstein inequality

\begin{equation}\label{bern}
||P_N w ||_{H^2(\R^2)}\leq ||u^m||_{H^2(\R^2)}+c(1+N^2)||v^m||_{L^2(\R^2)}.
\end{equation}

\noindent Therefore

\begin{equation}\label{bern2}
\limsup_m ||P_N w ||_{H^2(\R^2)}\leq \limsup_m ||u^m||_{H^2(\R^2)}.
\end{equation}

\noindent The sequence $P_Nw$ is bounded in $H^2(\R^2)$ and converges
towards $w$ in $H^1(\R^2)$. This completes the proof of the lemma.
\end{proof}

\section{The initial value problem in the energy space}

\subsection{Solutions for initial data in $\mathbb{E}_0=L^2(\R^2)\times L^2(\R^2) \times H^{-1}(\R^2)$.}

\noindent Consider the initial value problem for the system \eqref{e.1}-\eqref{nls}.
\underline{For local in time estimate} the dissipation parameters do not play a role
and we set $\alpha=\gamma=0$ for the sake of simplicity.
The Duhamel's form of the system reads

\begin{eqnarray}
\label{du1} \psi(t)=& S(t)\psi_0+i\int_0^t S(t-s)(\psi\phi-if )ds,\\
\label{du2} v(t)=& U(t)v(0)+\int_0^t U(t-s) i (Id-\Delta )^{-\frac12}(-Re \psi_{x_1} +g)ds,
\end{eqnarray}

\noindent where $v(t)=U(t)v_0$ is the solution of \eqref{nls} for $\gamma=0, g=0, \psi=0$
emanating from $v(0)$ at $t=0$.

\begin{lemma}\label{l2}
For any $(\psi_0,v_0)$ in $L^2(\R^2)\times L^2(\R^2)$ there exists a unique
solution of the equation in $C(0,T;L^2(\R^2))\cap L^4(0,T;L^4(\R^2))\times C(0,T;L^2(\R^2))$
for $T$ small enough depending on the $L^2$  norm of the initial data.
\end{lemma}

\begin{proof}

Set $Z=C(0,T;L^2(\R^2))\cap L^4(0,T;L^4(\R^2))$.
Set $G(\psi)(t)$ for the right hand side of \eqref{du1}, multiplied by
a cut-off function of time whose support is of size $T$.
Due to Strichartz estimate for the admissible pairs $(4,4)$ and $(\infty,2)$ we have

\begin{equation}\label{athens1}
||G(\psi)||_Z\leq c||\psi_0||_{L^2}+ ||\psi \phi||_{L^\frac43(0,T\times \R^2)}+ cT ||f||_{L^2(\R^2)}.
\end{equation}

\noindent Since by H\"older inequality

\begin{equation}\label{athens2}
||\psi \phi||_{L^\frac43(0,T\times \R^2)}\leq ||\phi||_{L^4(0,T\times \R^2} \sqrt T ||\psi||_{L^\infty(0,T;L^2)},
\end{equation}

\noindent for $T$ small enough we have that $G$ maps a ball of radius $2c||\psi_0||_{L^2}$ into itself, the $\psi$ term in the right hand
side of \eqref{athens2} being bounded by above by $c||\psi||_{L^\infty(0,T;L^2)}$ appealing
the linear equation \eqref{du2}. The proof that $G$ is a contraction is similar and then omitted.
\end{proof}

\begin{remark}
Actually this local solution is global in time. This will be proved in the absorbing set Section \ref{abs} below.
\end{remark}

\subsection{Solutions for initial data in $H^1(\R^2)\times H^1(\R^2) \times L^2(\R^2)$.}

For initial data in the energy space we just have to prove that the regularity
propagates along the flow. Computing the gradient in space of \eqref{du1} (multiplied by a cut-off function in time as above) we have,
using $$ \int_0^t S(t-s)\nabla fds= (S(t)-Id)(-\Delta)^{-1}\nabla f,$$
\noindent that

\begin{equation}\label{athens11}
||\nabla \psi||_{L^2}\leq ||\nabla \psi_0||_{L^2}+ c ||\nabla(\psi \phi)||_{L^\frac43(0,T\times \R^2}+ cT||f||_{L^2}.
\end{equation}

\noindent We have, due to H\"older inequality

$$ ||\phi \nabla \psi||_{L^\frac43(0,T\times \R^2}\leq T ||\nabla\psi||_{L^2} ||\phi||_{L^4}\leq T ||\nabla\psi||_{L^2} ||\phi||_{H^1}, $$

\noindent and

$$ ||\psi \nabla \phi||_{L^\frac43(0,T\times \R^2)}\leq ||\psi||_{L^4(0,T\times \R^2)}\sqrt T ||\nabla \phi||_{L^2}.$$

Since the linear equation \eqref{du2} propagates the regularity, we have the result.

\section{Bounded absorbing sets}
\label{abs}
\subsection{Absorbing set in $\mathbb{E}_0$}

We begin with an $L^2$ estimate for $\psi$. Consdering the scalar product
of the equation \eqref{e1.1} with $\psi$ we have

\begin{equation}\label{athens21}
\frac{d}{dt}||\psi||^2_{L^2(\R^2)}+2\alpha ||\psi||^2_{L^2(\R^2)}=2{\rm Re}\int_{\R^2} f \overline \psi dx \leq \alpha ||\psi||^2_{L^2(\R^2)}+ \frac{||f||^2_{L^2(\R^2)}}{\alpha}.
\end{equation}

\noindent Applying the Gronwall's lemma we have that the ball of radius $2\frac{||f||_{L^2(\R^2)}}{\alpha}$
is absorbing in $L^2$.

For the absorbing ball for $(\phi,\phi_t)$
we just have to prove the result for $(\phi-\phi^*,\phi_t)=(w,w_t)$ that is a solution to

\begin{equation}\label{athens22}
w_{tt}-\Delta w +w+\gamma w_t  =   -Re \psi_{x_1} .
\end{equation}

Introduce a new unknown $z$ such that

\begin{eqnarray}\label{athens23}
w_t+\frac\gamma 2 w=& z, \\
\label{athens235} z_{t}+\delta z-\Delta w +(1-\frac{\gamma^2}{4})w=&  -Re \psi_{x_1}.
\end{eqnarray}

\noindent Considering the scalar product of the first equation by $ w$
and of the second equation by $((1-\frac{\gamma^2}{4})Id-\Delta)^{-1}z$ we have

\begin{equation}\label{athens24}
\frac12\frac{d}{dt}(||w||^2_{L^2(\R^2)}+ ||z||^2_*)+\delta (||w||^2_{L^2(\R^2)}+ ||z||^2_*)\leq c||z||_*||\psi||_{L^2(\R^2)},
\end{equation}

\noindent where $||z||^2_*=\langle z,(1-\frac{\gamma^2}{4})(Id-\Delta)^{-1}z \rangle$ that is equivalent to the $H^{-1}$ norm of $z$.
Applying the Gronwall lemma to \eqref{athens24} concludes the proof of the absorbing ball in $\mathbb{E}_0$.

We can rephrase the last computation in (since the equation is linear)

\begin{lemma}\label{propabs}
For any $\sigma \leq 2$ (the limitation comes from the regularity of the stationary solution $\phi^*$).
If there exists an absorbing set for $\psi$ in $H^\sigma(\R^2)$
then there exists an absorbing set for $(\phi,\phi_t)$ in $H^\sigma(\R^2)\times H^{\sigma-1}(\R^2)$.
\end{lemma}

\subsection{Absorbing ball in the Energy space $\mathbb{E}_1$.}

We begin with an estimate on $\psi$.
Introduce the functional

\begin{equation}\label{athens31}
H=\int_{\R^2}( |\nabla \psi|^2 +\phi |\psi|^2+ 2{\rm Re}(f\overline \psi))dx,
\end{equation}

\noindent and considering the scalar product of \eqref{e1.1} with $i(\psi_t+\alpha \psi)$ leads to

\begin{equation}\label{athens32}
\frac{d}{dt} H +\alpha H= \int_{\R^2} \phi_t |\psi|^2dx + {\rm Re}\int_{\R^2}(f\overline \psi)dx.
\end{equation}

\noindent We then have

\begin{equation}\label{athens33}
\frac{d}{dt} H +\alpha H \leq ||\phi_t||_{L^2(\R^2)}||\psi||^2_{L^4(\R^2)}+ K\leq \varepsilon ||\phi_t||_{L^2(\R^2)}^2+ \frac{||\psi||^4_{L^4(\R^2)}}{\varepsilon}+K .
\end{equation}

\noindent Let us observe that moreover on the $\mathbb{E}_0$ absorbing set, by H\"older and Gagliardo-Nirenberg inequalities

\begin{equation}\label{Hcoercivity}
H\geq  ||\nabla \psi||^2_{L^2(\R^2)}-||\phi||_{L^2(\R^2)}||\psi||_{L^2(\R^2)}||\nabla \psi||_{L^2(\R^2)}-\frac K 2\geq      \frac12 ||\nabla \psi||^2_{L^2(\R^2)}-K.
\end{equation}

We set $N=<w,(1-\frac{\gamma^2}{4})(Id-\Delta) w>+||z||^2_{L^2(\R^2)}$.

\noindent Considering the scalar product of
\eqref{athens23} by $(1-\frac{\gamma^2}{4})(Id-\Delta) w$
and of ({\color{red} 21}) by $z$, we have, summing the resulting equations

\begin{equation}\label{athens34}
\frac{d}{dt} N +\frac\gamma 2 N\leq ||z||_{L^2(\R^2)}||\nabla \psi||_{L^2(\R^2)}\leq \frac\gamma 4 N+ \frac C \gamma H.
\end{equation}

\noindent Adding \eqref{athens33} and \eqref{athens34} yields, for $\varepsilon$ small enough, and setting $2\delta=\min(\alpha,\frac\gamma 4)$,

\begin{equation}\label{athens35}
\frac{d}{dt}(H+\frac{\alpha\gamma}{2C} N)+\delta (H+\frac{\alpha\gamma}{2C} N) \leq  K (1+ ||\psi||^4_{L^4(\R^2)}).
\end{equation}

We then have, for $J=H+\frac{\alpha\gamma}{2C} N$ that is the square of a norm in $\mathbb{E}_1$,

\begin{equation}\label{athens36}
J(t)e^{\delta t}\leq J(0)+K \int_0^t e^{\delta s} (1+ ||\psi||^4_{L^4(\R^2)})ds.
\end{equation}

\noindent Assume now that the trajectory is in the $\mathbb{E}_0$ absorbing ball for $t\geq 0$
(up to a time translation), fixing $T$ as in the proof of Lemma 2.2., we have for $NT\simeq t$

\begin{equation}\label{athens37}
\int_0^t e^{\delta s} ||\psi||^4_{L^4(\R^2)}ds\leq c \sum_{k=0}^N e^{\delta kT}\int_{kT}^{kT+T} ||\psi||^4_{L^4(\R^2)}ds\leq K e^{\delta t}.
\end{equation}

\noindent Then the proof of the absorbing ball in $\mathbb{E}_1$ is complete.

\section{Proof of the existence of a global attractor}
\label{proofa}

We begin with an argument due to \cite{LuWang}.

\begin{lemma}\label{luwang}
For any $\varepsilon>0$ there exists $R>0$ such that for $t$ large enough uniformly in the $\mathbb{E}_1$
absorbing set $\int_{|x|>R} |\psi|^2dx \leq \varepsilon$.
\end{lemma}

\begin{proof}

Consider a non negative smooth function $\theta$ such that $\theta(x)=1$ if $|x|\leq 1$ and $\theta(x)=0$ if $|x|>2$.
Set $\theta_R(x)=\theta(\frac{| x|}{ R})$. We define the seminorm representing the tail of $\psi$ as

\begin{equation}\label{seminorm}
||\psi||_{\theta}=(\int_{\R^2} (1-\theta_R(x))|\psi(x)|^2dx)^\frac12.
\end{equation}

Considering the scalar product of \eqref{e1.1} with $(1-\theta_R)\psi$ leads to

\begin{equation}\label{tai1}
\frac12 \frac{d}{dt} ||\psi||_{\theta}^2+\alpha ||\psi||_{\theta}^2= {\rm Im}\int_{\R^2} \Delta \psi(1-\theta_R)\overline \psi +{\rm Re}\int_{\R^2} f(1-\theta_R)\overline \psi.
\end{equation}

\noindent The second term of the right hand side of \eqref{tai1} is bounded by above, thanks to Cauchy-Schwarz inequality, by
$||\psi||_{\theta} (\int_{|x|>R}|f|^2)^\frac12 $. Integrating by parts the first term in the right hand side of \eqref{tai1} we have

\begin{equation}\label{ta2}
| {\rm Im}\int_{\R^2} \Delta \psi(1-\theta_R)\overline \psi|=|R^{-1}\int_{\R^2} \theta'(\frac{|x|}{ R}) {\rm Im}(\frac{x}{|x|}.\nabla \psi(x)\overline \psi(x))dx|\leq R^{-1}||\theta'||_{L^\infty}||\psi||^2_{H^1(\R^2)}.
\end{equation}

\noindent Since $\psi$ is included in the global attractor and then remains bounded in $H^1(\R^2)$ we have

\begin{equation}\label{ta3}
\frac{d}{dt} ||\psi||_{\theta}^2+\alpha ||\psi||_{\theta}^2\leq \frac 1 \alpha \int_{|x|>R}|f|^2+ \frac K R.
\end{equation}

\noindent We conclude by the Gronwall Lemma.
\end{proof}

As a mere consequence of this tail estimate

\begin{corollary}\label{asycom}
Any set of the first component of a trajectory $\psi(t)$
that is bounded in $H^1(\R^2)$ and that satisfies the estimate
in Lemma \ref{luwang} is asymptotically compact in $L^2(\R^2)$.
\end{corollary}

We now complete the proof of the existence of a global attractor.
The flow map $(\psi_0,v_0)\mapsto (\psi,v)$ is continuous
for the weak topology of $\mathbb{E}_1$.
To pass to the limit in \eqref{e1.1} the difficulty is to pass to the limit
in the non-linear term $\phi\psi$; due to corollary \eqref{asycom}
the convergence of the $\psi$ term  holds strongly in $L^2$
and by interpolation in any $L^p$ space, $p<+\infty$. Then
we can pass to the limit for this nonlinear term. The $v$
part of the equation is linear and then pass to the limit.

Following \cite{Ghidaglia}, the restriction of the flow
map to the absorbing set in $\mathbb{E}_1$ is a compact
mapping for the weak topology, and we then have
a weak global attractor $\mathcal{A}$  in $\mathbb{E}_1$.
Due to the famous J. Ball's Lemma we can prove that
$\mathcal{A}$ is in fact a global attractor in the energy space
for the semigroup.

\section{Regularity of the global attractor}

We begin with a $H^2$ estimate for $\psi$.

\subsection{Extra regularity on $\psi$}
\label{firststep}

We prove here that the first component $\psi(t)$ of
a  complete bounded trajectory $(\psi(t),\phi(t),\phi_t(t))$ in the energy space
is in fact included in a bounded set of $H^2(\R^2)$.

Let $m$ be in $\mathbb{N}$. We split $\psi(t)=u^m(t)+v^m(t)$ where $u^m$ is solution for $t\geq -m$
to the equation \eqref{mp1} below with non standard initial condition at $t=-m$
that reads $u^m(-m)=u_m$ where

\begin{eqnarray}
\label{mp1} u^m_t +i\Delta u^m  +\alpha u^m &= & i\phi u^m+f, \\
\label{mp1.0}\alpha u_m+i \Delta u_m &= &i\phi(-m)u_m+f.
\end{eqnarray}

\noindent For a pair $(\phi(-m),f)$ there exists a unique
$u_m$ solution of the elliptic equation \eqref{mp1.0} thanks to Lax-Milgram theorem.
Then the remaining term $v^m=\psi-u^m$ solves

\begin{equation}\label{mp2}
 v^m_t +i\Delta v^m  +\alpha v^m =  i\phi v^m ,
\end{equation}

\noindent with initial data $v^m(-m)=v_m=\psi(-m)-u_m$ at $t=-m$.
We now establish that $v^m$ vanishes when $m$ goes to $+\infty$.
Considering the scalar product of \eqref{mp2} with $v$ and integrating in $t$ between $-m$ and $0$ leads to

 \begin{equation}\label{mp3}
 ||v^m(0)||^2_{L^2(\R^2)}\leq ||v^m(-m)||^2_{L^2(\R)^2}\exp(-2\alpha m)\leq K \exp(-2\alpha m).
 \end{equation}

 We now prove that $u^m$ belongs to a more regular space than $H^1(\R^2)$.
 Differentiating \eqref{mp1} with respect to $t$ and setting $w^m=u^m_t$ leads to, thanks to the initial condition \eqref{mp1.0}

 \begin{eqnarray}
\label{mp4} w^m_t +i\Delta w^m  +\alpha w^m &= & i\phi w^m+i \phi_t u^m, \\
\label{mp4.0}w^m(-m)&= &0.
\end{eqnarray}

Considering the scalar product of \eqref{mp4} by $w^m$ leads to

 \begin{equation}\label{mp5}
 \frac12 \frac{d}{dt} ||w^m||^2_{L^2(\R^2)}+\alpha ||w^m||^2_{L^2(\R^2)}=-{\rm Im}\int_{\R^2} \phi_t u^m \overline{w}^m.
 \end{equation}

\noindent Since the right hand side of \eqref{mp5} is bounded by

$$||w^m||_{L^2(\R^2)} ||\phi_t||_{L^2(\R^2)} ||u^m||_{L^\infty(\R^2)}\leq K ||w^m||_{L^2(\R)}||u^m||_{L^\infty(\R^2)}.$$

\noindent Thanks to Agmon's inequality and since $\Delta u^m=iu^m_t+$ lower order terms
we have

$$ ||u^m||_{L^\infty(\R^2)}\leq c||u^m||_{L^2(\R^2)}^\frac12 ||\Delta u^m||_{L^2(\R^2)}^\frac12 \leq K(1+ ||w^m||_{L^2(\R^2)}^\frac12).$$

\noindent Therefore applying Gronwall lemma leads to $||w^m(0)||_{L^2(\R^2)}\leq K .$

We know that $\psi(0)=u^m(0)+v^m(0)$ where $u^m(0)$ remains in a bounded set of $H^2(\R^2)$
and $v^m(0)$ converges to $0$ in $L^2(\R^2)$. Then we have a uniform bound in $H^2(\R)$ for $\psi(0)$.
Using the same argument replacing $t=0$
by any $t$ yields that the first component $\psi(t)$ of a complete trajectory remains in a bounded set of $H^2(\R^2)$.

\subsection{Extra regularity on $\phi,\phi_t$}

The $H^2(\R^2)\times H^1(\R^2)$ bound on $(\phi,\phi_t)$
is a consequence of Lemma \ref{lemma}, setting
$\sigma=1$, $u=\phi-\phi^*$ and $h={\rm Re}(\psi_{x_1}^*-\psi_{x_1})$.

\subsection{Extra regularity on $\psi_t$}
\label{thirdstep}

Consider a complete trajectory in the global attractor.
The equation \eqref{e1.1} and the two previous subsections yield that $\psi_t$
remains in a bounded set of $L^2(\R^2)$. In fact this result can be improved

\begin{proposition}\label{inexpected}
Consider a complete trajectory in the global attractor. Then $\psi_t$ remains bounded in $H^1(\R^2)$.
\end{proposition}

\begin{proof}
We use the splitting $\psi=u^m+v^m$ introduced in Subsection \ref{firststep}.
The part $v^m$ is the vanishing part when $m$ diverges towards $+\infty$;
 the derivative $v^m_t$ vanishes also thanks to \eqref{mp2}.
We focus on $w^m=u^m_t$. Considering the scalar product of \eqref{mp4} with $-i(w^m_t+\alpha w^m)$ leads
to, for $t\geq -m$

\begin{equation}\label{in1}
\frac 1 2 \frac{d}{dt} J(w^m)+\alpha J(w^m)=\frac 32 \int_{\R^2} \phi_{t} |w^m|^2 +{\rm Re} \int_{\R^2} \phi_{tt} u^m \overline{ w}^m+
\alpha {\rm Re} \int_{\R^2} \phi_t u^m \overline{w}^m ,
\end{equation}

\noindent where $J(w^m)=||\nabla w^m||^2_{L^2(\R^2)}+\int_{\R^2} \phi |w^m|^2 +2{\rm Re}\int_{\R^2} \phi_t u^m \overline{w}^m.$
Due to the previous estimates on the global attractor

\begin{eqnarray}
J(w^m)&\geq& ||\nabla w^m ||^2_{L^2(\R^2)}- (||\phi||_{L^\infty}+2||\phi_t||_{L^4(\R^2)}||u^m||_{L^4(\R^2)}) ||w^m||^2_{L^2(\R^2)}\\
\label{in2} &\geq& ||\nabla w^m||^2_{L^2(\R^2)}-  K.
\end{eqnarray}

\noindent We now handle the right hand side of \eqref{in1}. We have, thanks to Gagliardo-Nirenberg inequality

\begin{equation}\label{in3}
\frac 12 \int_\R \phi_{t} |w^m|^2 \leq ||\phi_{t}||_{L^2(\R^2)}||w^m||^2_{L^4(\R^2)}\leq K ||\nabla w^m ||_{L^2(\R^2)}.
\end{equation}

\noindent Analogously

\begin{equation}\label{in5}
|\int_\R \phi_{tt} u^m \overline{ w}^m|\leq ||\phi_{tt}||_{L^2(\R^2)}||u^m||_{L^\infty(\R^2)}||w^m||_{L^2(\R^2)}\leq K,
\end{equation}

\noindent and

\begin{equation}\label{in7}
|\int_\R \phi_t u^m \overline{w}^m| \leq ||\phi_{t}||_{L^2(\R^2)}||u^m||_{L^\infty(\R^2)}||w^m||_{L^2(\R^2)}\leq K.
\end{equation}

\noindent Gathering these inequalities the Gronwall Lemma completes the proof of the Proposition.

\end{proof}

\begin{corollary}\label{extra2}
Consider a complete trajectory in the global attractor. Then $\phi_t$ remains bounded in $H^2(\R^2)$.
\end{corollary}

\begin{proof}
Differentiate \eqref{e1.2} with respect to $t$ and apply Lemma \ref{lemma} with $\sigma=1$ to
the equation

\begin{equation}\label{extra3}
(\phi_t)_{tt}+\gamma (\phi_t)_t+\phi_t-\Delta \phi_t=-{\rm Re}(\psi_{tx_1}).
\end{equation}

\end{proof}

\subsection{Compactness of the global attractor in $H^2(\R^2)\times H^2(\R^2) \times H^2(\R^2)$}

We now state

\begin{lemma}\label{pan1}
The set $\psi=\psi(0)$ such that $(\psi(t),\phi(t), \phi_t(t))$ belongs to $\mathcal{A}$
is a compact set of $L^2(\R^2)$.
\end{lemma}

\begin{proof}

We use here a splitting of the solution as in the proof of the existence
of a global attractor in \cite{Laurencot}.
For $\varepsilon>0$ consider a compactly supported function $f_\varepsilon$ in $H^1(\R^2)$ such that
$||f-f_\varepsilon||_{L^2(\R)}\leq \varepsilon$. Consider now the splitting $\psi=u^{m,\varepsilon}+v^{m,\varepsilon}$
where for $t\geq -m$,

\begin{eqnarray}
\label{mf1} u^{m,\varepsilon}_t +i\Delta u^{m,\varepsilon}  +\alpha u^{m,\varepsilon} &= & i\phi u^{m,\varepsilon}+f_\varepsilon, \\
\label{mf1.0}u^{m,\varepsilon}(-m)&= &0.
\end{eqnarray}

\noindent and

\begin{eqnarray}\label{mf2}
 v^{m,\varepsilon}_t +i\Delta v^{m,\varepsilon}  +\alpha v^{m,\varepsilon} &= & i\phi v^{m,\varepsilon}+f-f_\varepsilon, \\
\label{mf2.0} v^{m,\varepsilon}(-m)&= &\psi(-m).
\end{eqnarray}

\noindent Considering the scalar product of \eqref{mf2} with $ v^{m,\varepsilon} $ yields

\begin{equation}\label{mf3}
\frac12 \frac{d}{dt} ||v^{m,\varepsilon}||^2_{L^2(\R^2)}+\alpha  ||v^{m,\varepsilon}||^2_{L^2(\R^2)}\leq ||f-f_\varepsilon||_{L^2(\R^2)} ||v^{m,\varepsilon}||_{L^2(\R^2)},
\end{equation}

\noindent and then by Gronwall Lemma

\begin{equation}\label{mf4}
||v^{m,\varepsilon}(0)||_{L^2(\R^2)}\leq \exp(-\alpha m) ||\psi(-m)||_{L^2(\R^2)}+\frac \varepsilon \alpha.
\end{equation}

\noindent Similarly to the computations in Subsection \ref{firststep} we can prove that $v^{m,\varepsilon}$
remains in a bounded set of $H^1(\R^2)$. Then $u^{m,\varepsilon}$ remains in a bounded set of $H^1(\R^2)$ as well.

\noindent Considering the scalar product of \eqref{mf1} with $|x|^2 u^{m,\varepsilon}$ yields

\begin{eqnarray}\label{mf3}
\frac12 \frac{d}{dt} ||xu^{m,\varepsilon}||^2_{L^2(\R^2)}+\alpha  ||xu^{m,\varepsilon}||^2_{L^2(\R^2)}&\leq & |{\rm Im}\int_{\R^2} 2x.\nabla u^{m,\varepsilon} \overline{u}|
+|{\rm Re}\int_{\R^2} |x|^2f_\varepsilon \overline{u}^{m,\varepsilon}|\\  &\leq & \frac \alpha 2 ||xu^{m,\varepsilon}||^2_{L^2(\R^2)}+ K+ ||xf_\varepsilon||_{L^2(\R^2)}^2.
\end{eqnarray}

\noindent Thanks to Gronwall lemma we have that $xu^{m,\varepsilon}(0)$ remains in a bounded set of $L^2(\R^2)$,
the bound depending on $\varepsilon$ but not on $m$. Therefore letting $m$ diverge to $+\infty$
we have that $\psi(0)$ is captive in $K_\varepsilon+B_E(0,\frac \varepsilon \alpha)$ where $K_\varepsilon$
is closed and bounded in $H^1(\R^2) \cap L^2(|x|^2dx)$ and then compact in $L^2(\R^2)$.
\end{proof}

By interpolation with the $H^2(\R^2)$ bound  we have

\begin{corollary}
Set $\sigma<2$. The set $\psi=\psi(0)$ such that $(\psi(t),\phi(t), \phi_t(t))$ belongs to $\mathcal{A}$
is a compact set of $H^\sigma(\R^2)$.
\end{corollary}

We now propagate the compactness to the $(\phi,\phi_t)$ part of
the attractor.

\begin{lemma}\label{pan2}
For $\sigma<2$, the set $(\phi,\phi_t)=(\phi(0),\phi_t(0))$ such that $(\psi(t),\phi(t), \phi_t(t))$ belongs to $\mathcal{A}$
is a compact set of $H^\sigma(\R^2)\times H^{\sigma}(\R^2)$.
\end{lemma}

\begin{proof}
We use the representation formula \eqref{rep}

\begin{equation}\label{rep2}
\phi(0)+i(Id-\Delta)^{-\frac12}\phi_t(0)=i\int_{-\infty}^0 U_{\gamma}(-s) (Id-\Delta)^{-\frac12}(-{\rm Re}\psi_{x_1}(s)+g)ds.
\end{equation}

Consider $\sigma\geq 1$. The map $\psi(0) \mapsto \psi(s)$ is continuous from $H^{\sigma}(\R^2)$
to $C([-T, 0]; H^\sigma(\R^2))$. This is a byproduct of solving the initial value problem.
Fix now $T$ large enough such that

$$ ||\int_{-\infty}^{-T} U_{\gamma}(-s) (Id-\Delta)^{-\frac12}(-{\rm Re}\psi_{x_1}(s)+g)ds||_{H^\sigma(\R^2)}\leq \varepsilon.$$

\noindent This part is then included in a small ball.
\noindent Besides, since the map $$F: \psi(0)\mapsto \int_{-T}^0 U_{\gamma}(-s) (Id-\Delta)^{-\frac12}(-{\rm Re}\psi_{x_1}(s)+g)ds$$
\noindent is continuous, then $F$ maps a compact set into a compact set.
Thanks to Lemma \ref{compact} the pair $(\phi(0),\phi_t(0))$ remains in compact set of $H^\sigma(\R^2)\times H^{\sigma-1}(\R^2)$.
For the $H^\sigma(\R^2)$ compactness we interpolate with the $H^2(\R^2)$ bound for $\phi_t$.
\end{proof}

To reach the compactness of $\psi$ for $\sigma=2$, we prove that $w=\psi_t$
is in a compact set of $L^2(\R^2)$. Differentiating \eqref{e1.1} with respect to
$t$ and considering the scalar product of the resulting equation with $w$ yields

\begin{equation}\label{Ball}
||w(t)||^2_{L^2(\R)}=||w(0)||^2_{L^2(\R)}e^{-2\alpha t}-\int_0^t e^{2\alpha(s- t)}(\int_\R {\rm Im} \phi_t(s)\psi(s)\overline{w}(s)dx) ds.
\end{equation}

The right hand side is weakly continuous in $L^2(\R^2)$ since both $\phi_t$ and $\psi$
are compact in $H^1(\R^2)$. Then  the framework of the J. Ball's argument applies
\cite{Ball}, \cite{mrw}, \cite{Wang} and we have the compactness of the $\psi$ part of the attractor in $H^2(\R^2)$.
For the $(\phi,\phi_t)$ part we propagate the compactness using Lemma \ref{pan2}.
This completes the proof of Theorem 1.1.

\section{The attractor has finite dimension}

In this section we prove that the global attractor has finite
fractal dimension.
 We follow here the Douady-Oesterl\'e \cite{DO} method
as described in \cite{Ghidaglia}, \cite{Temam} of contraction
of infinitesimal volumes. The global attractor is a compact subset of $ H^2(\R^2)\times  H^2(\R^2)\times  H^2(\R^2)$.
Since $H^2(\R^2)$ is a Banach algebra we can differentiate the semi-group
defined by the flow map of solutions in $\mathcal{A}$. For alternate proofs that do not use the differentiability of the semigroup
we refer to \cite{zelik}, \cite{chueshov}.

\subsection{A general result}

Let $X$ be a separable Hilbert space.
Set respectively $(u,v)_X$ and $||u||_X$ for respectively the scalar product and norm in $X$.
For a collection $(u_1,...u_m)$ of $m$ vectors in $X$
set $G(u_1,....u_m)$ for the Gram determinant of this collection of vectors,
i.e. the determinant of the matrix $g$ whose entries are $g_{kl}=(u_k,u_l)_X$.
This quantity is the volume of the parallelepiped whose vertices are $(u_1,...u_m)$.

\begin{theorem}\label{chulas}
Consider $(h_1,....h_m)$ an orthonormal collection of $m$ vectors in $X$.
Let $\mathcal{A}$ be a compact set in $X$ that is invariant by a semigroup $S(t)$ acting on $X$.
For $a$ in $\mathcal{A}$, consider $S(t)a$ the complete trajectory emanating from $a$
at $t=0$ and $u_j(t)=DS(t)(a)h_j$ the derivative of the semigroup emanating from $h_j$
at $t=0$. If the $m-$dimensional volume $G(u_1(t),....u_m(t))$ decays to $0$ when
$t$ goes to $+\infty$ then $\dim \mathcal{A}\leq m$.
\end{theorem}

We use the following corollary to conclude.
Let $Y$ be another space such that the embedding $X \subset Y$ is compact.
Set $(u,v)_Y$ and $||u||_Y$ for the scalar product and norm in $Y$.

\begin{corollary}
Let $\delta$ be a positive number
If we have the following contraction estimate, up to a compact term
\begin{equation}\label{boa}
\frac{d}{dt} ||u||^2_X +\delta ||u||^2_X \leq K ||u||^2_Y,
\end{equation}
\noindent then the dimension of $\mathcal{A}$ is finite.
\end{corollary}

\begin{proof}
Since $X\subset Y$ is compact there exists a self-ajoint compact operator $T$
such that $(u,v)_Y=(u,Ty)_X$ for any vectors $u,v$ in $X$.
Moreover there exists a countable sequence of eigenvalues for $T$
that satisfies $\lambda_j\geq \lambda_{j+1}\rightarrow 0$.
Introduce $G(t)=G(u_1(t),..u_m(t))$. Then we infer from \eqref{boa} by the Courant min-max formula
(see \cite{Ghidaglia}, \cite{Temam} and the references therein)
\begin{equation}
\frac{d G}{dt} + 2\delta m G \leq K \left(\sum_{j=1}^m \lambda_j\right) G.
\end{equation}

\noindent Since $\frac 1 m \sum_{j=1}^m \lambda_j $ converges to $0$ when $m$
diverges to $+\infty$, then the volume $G(t)$ decays to $0$ when $t$ diverges to $+\infty$
if $m$ is such that $\frac 1 m \sum_{j=1}^m \lambda_j  < \frac{2\delta}{ K}$.
\end{proof}

\subsection{Proof of Theorem \ref{th2}}

We follow now ideas introduced by B. Alouini (see \cite{Alouini} and the references therein).
The strategy reads as follows. Thanks to Rellich Theorem, the embedding $H^1(\R^2) \subset L^2([-R,R]^2)$ is compact.
The first step in Section \ref{firststepbeyond} below is then to obtain uniform bounds on the tails of the solution, that
are the semi-norms of a solution outside a given interval $[-R,R]$.
The second step in Section \ref{secondstepbeyond} consists in proving the key estimate \eqref{boa}
with $X$ the energy space and $Y=L^2([-R,R]^2)\times \{0\} \times \{0\}$,
using that the map $(\psi,\phi,\phi_t) \mapsto (\psi,0,0)$ is compact.

\subsubsection{Estimates on the tails of a trajectory}
\label{firststepbeyond}

Consider the cut-off function $\theta_R(x)$ as in Section \ref{proofa}.
We consider semi-norms defined in \eqref{seminorm} above.

\begin{lemma}\label{tail1}
For any $\varepsilon>0$ there exists a $R$ large enough depending on the bounds on the attractor such that
\begin{equation}\label{ta0}
||\psi_t||_\theta \leq \varepsilon.
\end{equation}
\end{lemma}

\begin{proof}

Set $w=\psi_t$ that is solution to

\begin{equation}\label{orkyn}
w_t+\alpha w +i\Delta w=i\phi w+i \phi_t \psi,
\end{equation}

\noindent that is defined for any $t$ in $\R$ and bounded in $H^1(\R^2)$.
Considering now the scalar product of \eqref{orkyn} with $(1-\theta_R) w$
leads to

\begin{equation}\label{ta4}
\frac12 \frac{d}{dt} ||w||_{\theta}^2+\alpha ||w||_{\theta}^2= {\rm Im}\int_{\R^2} \Delta w(1-\theta_R)\overline w -{\rm Im}\int_{\R^2} \phi_t (1-\theta_R)\overline w \psi.
\end{equation}

\noindent The first term in the right hand side of \eqref{ta4} is bounded by above, proceeding as in \eqref{ta2}, by $R^{-1}||w||^2_{H^1(\R^2)}$.
The second term is bounded by $||\phi_t||_{L^\infty(\R^2)}||\psi||_\theta||w||_\theta$. Therefore

\begin{equation}\label{ta5}
\frac{d}{dt} ||w||_{\theta}^2+\alpha ||w||_{\theta}^2\leq R^{-1}||w||^2_{H^1(\R^2)}+\frac 1 \alpha ||\phi_t||_{L^\infty(\R^2)}^2||\psi||_\theta^2.
\end{equation}

\noindent Using Proposition \ref{inexpected}, the Lemma \ref{luwang}, and
integrating between $-\infty$ and $0$ completes the proof.
\end{proof}

We proceed similarly for the $(\phi,\phi_t)$.

\begin{lemma}\label{tail2}
For any $\varepsilon>0$ there exists a $R$ large enough depending on the bounds on the attractor such that
\begin{equation}\label{ta00}
||\phi||_\theta+||\phi_x||_\theta+ ||\phi_t||_\theta \leq \varepsilon.
\end{equation}
\end{lemma}

\begin{proof}

Let us recall that $0<\gamma<<1$.
Following an idea introduced in \cite{ghidagliatemam}, the equation \eqref{e1.2} reads also

\begin{eqnarray}
 \label{e1.21} \phi_t +\frac \gamma 2 \phi &= & y, \\
 \label{e1.22}
 y_{t}+ \frac \gamma 2 y -\Delta \phi +(1-\frac{\gamma^2}{4})\phi &= &  -Re \psi_{x_1} +g.
\end{eqnarray}

\noindent We consider the scalar product of the first equation by $(1-\theta_R)( -\Delta \phi +(1-\frac{\gamma^2}{4})\phi)$
and the scalar product of the second equation by $(1-\theta_R)y$. Summing the resulting equations and
setting $I= ||y||_{\theta}^2+||\phi_x||_\theta^2+(1-\frac{\gamma^2}{4})||\phi||^2_\theta$, we obtain

\begin{equation}\label{ta11}\begin{split}
\frac12 (\frac{dI}{dt}   +  \gamma I) \leq
-{\rm Re}\int_{\R^2} \psi_{x_1} (1-\theta_R) y+\frac 1 R \int_{\R^2} \theta'(\frac x R) y \frac{x}{|x|}.\nabla \phi+ \int_{\R^2} g (1-\theta_R) y.
\end{split}\end{equation}

\noindent We have then

\begin{equation}\label{ta12}
\frac{dI}{dt}   + \frac \gamma 2 I\leq \frac{8}{\gamma}(||\psi_{x_1}||^2_{|x|>R}+||g||^2_{|x|>R})+\frac K R.
\end{equation}

\noindent Due to Lemma \ref{tail1} the conclusion follows promptly.

\end{proof}

\subsection{Checking the assumptions of Theorem \ref{chulas}}
\label{secondstepbeyond}

Consider $(\psi,\phi,\phi_t)$ a complete trajectory that belongs to the global attractor.
Introduce $(u,v,v_t)$ that are the differential of $S(t)$ along the trajectory.
Then $(u,v,v_t)$ solves the system, setting $z=v_t+\frac \gamma 2 v$

\begin{eqnarray}
 \label{f1.1} u_t +i\Delta u  +\alpha u  &=&  i\phi u +i \psi v, \\
 \label{f1.2} v_t+\frac \gamma 2 v & = & z, \\
 \label{f1.3} z_t+\frac \gamma 2 z -\Delta v + (1-\frac{\gamma^2}{4}) v  &=&  -Re u_{x_1}.
\end{eqnarray}

We begin with an $L^2$ estimate. Considering the scalar product of \eqref{f1.1} with $u$ yields,
using $1= (1-\theta_R)+\theta_R$,

\begin{eqnarray}\label{che1}
\frac12 \frac{d}{dt}||u||^2_{L^2(\R^2)}+\alpha ||u||^2_{L^2(\R^2)}&=&-{\rm Im}\int_{\R^2} (1-\theta_R) \psi v \overline u- {\rm Im}\int_{\R^2} \theta_R \psi v \overline u
\\ &\leq & K||v||_{L^2(\R^2)}(||u||_{L^2([-R,R]^2)} +||\psi||_\theta ||u||_{L^2(\R^2)})\\
&\leq & K||v||_{L^2(\R^2)}(||u||_{L^2([-R,R]^2)} +\varepsilon ||u||_{L^2(\R^2)}).
\end{eqnarray}

\noindent where $\varepsilon>0$ is defined in Lemmas \ref{luwang}, \ref{tail1}, \ref{tail2}. Introduce now

$$H(u)=||\nabla u||^2_{L^2(\R^2)}+\int_{\R^2} \phi |u|^2+2 {\rm Re}\int_{\R^2} v \psi \overline u.$$

\noindent Considering the scalar product of \eqref{f1.1} by $-i(u_t+\alpha u)$ yields

\begin{equation}\label{che2}
\frac12 \frac{d}{dt}H(u)+\alpha H(u)=\frac12 \int_{\R^2} \phi_t |u|^2 +\alpha {\rm Re}\int_{\R^2} v \psi \overline u+{\rm Re}\int_{\R^2} v_t \psi \overline u
+{\rm Re}\int_{\R^2} v \psi_t \overline u.
\end{equation}

\noindent To bound by above the right hand side of \eqref{che2} we use the $1=(1-\theta_R)+\theta_R$ trick.
For the first term we have

\begin{equation}\label{che3}
\int_{\R^2} \phi_t |u|^2 \leq ||\phi_t||_\theta ||u||^2_{H^1(\R^2)} + ||\phi_t||_{L^\infty(\R^2)} ||u||^2_X.
\end{equation}

\noindent We also have, for the third and fourth terms

\begin{equation}\label{che5}
{\rm Re}\int_{\R^2} v_t \psi \overline u \leq K||v_t||_{H^1(\R^2)}( \varepsilon ||u||_{L^2(\R^2)}+ ||u||_{X}),        ,
\end{equation}

\noindent and

\begin{equation}\label{che6}
 {\rm Re}\int_{\R^2} v \psi_t \overline u \leq  K||v||_{L^2(\R^2)}(\varepsilon ||u||_{H^1(\R^2)}+ ||u||_{X} ).
\end{equation}

Gathering these inequalities leads to, for $\kappa>0$ and $H_\kappa(u)=H(u)+\kappa ||u||^2_{L^2(\R^2}$,

\begin{equation}\label{gather}
\frac12 \frac{d}{dt}H_\kappa(u)+\alpha H_\kappa(u)\leq K(||v||_{L^2(\R^2)}+||u||_{H^1(\R^2)}+||v_t||_{L^2(\R^2)} )(||u||_{X}+\varepsilon ||u||_{H^1(\R^2)} ).
\end{equation}

Introduce now

\begin{equation}\label{N}
N= ||\nabla v||^2_{L^2(\R^2)}+(1-\frac{\gamma^2}{4})||v||_{L^2(\R^2)}^2+||z||^2_{L^2(\R^2)}.
\end{equation}

\noindent We have, considering the scalar product of \eqref{f1.2} with $-\Delta v + (1-\frac{\gamma^2}{4})v$ and the scalar
 product of \eqref{f1.3} with $z$ and summing the resulting equations

\begin{equation}\label{che9}\begin{split}
\frac12\frac{d  N}{dt} + \frac \gamma 2   N
  \leq -{\rm Re}\int_\R u_{x_1} z
\leq \frac \gamma 4 ||z||^2_{L^2(\R^2)}+ \frac1 \gamma ||\nabla u||_{L^2(\R^2)}^2.
\end{split}\end{equation}

Considering $\beta>0$ another parameter we have

\begin{equation}\label{dim1}\begin{split}
\frac{d}{dt} (H_\kappa+\beta N)+2\alpha H_\kappa +\beta \frac \gamma 2 N
\leq \frac \beta \gamma ||\nabla u||_{L^2(\R^2)}^2+\\ K(||v||_{L^2(\R^2)}+||u||_{H^1(\R^2)}+||v_t||_{L^2(\R^2)} )(||u||_{X}+\varepsilon ||u||_{H^1(\R^2)}).
\end{split}\end{equation}

We now chose $\beta$ small enough to ensure $3\alpha > \frac \beta \gamma$.
Then the term $\frac \beta \gamma ||\nabla u||_{L^2(\R^2)}^2$ can be absorbed in the left hand side.
Besides

\begin{equation}\label{dim2}
H(u)\geq \frac12 ||\nabla u||_{L^2(\R^2)}^2-K_1||u||^2_{L^2(\R^2)}-K_2||v||_{L^2(\R^2)}||u||_{L^2(\R^2)}.
\end{equation}

We use $K_2||v||_{L^2(\R^2)}||u||_{L^2(\R^2)}\leq \frac\beta 8 ||v||_{L^2(\R^2)}^2+ \frac{K_3}{\beta}||u||^2_{L^2}$,
then we chose $\kappa$ large enough such that $\kappa \geq 2K_1 +2 \frac{K_3}{\beta} $.
Then the quantity $H_\kappa+\beta N$ is (the square) of a norm in $\mathbb{E}_1$, say $||(u,v,z)||_Y^2$.
 Equation \eqref{dim1} reads

\begin{equation}\label{dim3}
\frac{d}{dt} (H_\kappa+\beta N)+\min(\alpha, \frac \gamma 2)(H_\kappa+\beta N)\leq \frac
K\varepsilon (||v||_{L^2(\R^2)}+||u||_{H^1(\R^2)}+||v_t||_{L^2(\R^2)} )^2 +\frac{\tilde K}{\varepsilon}||u||^2_X.
\end{equation}

\noindent We now chose $\varepsilon>0$ such that the first term in the right hand side
of \eqref{dim3} can be bounded by above by $\frac12 \min(\alpha, \frac \gamma 2)(H_\kappa+\beta N)$.
This completes the proof of \eqref{boa}.

\section*{Acknowledgments} The first author acknowledges the hospitality of the University of West Attica
and the support of the CDP C2EMPI, together with the French State under the France-2030 programme, the University of Lille, the Initiative of Excellence of the University of Lille, the European Metropolis of Lille for their funding and support of the R-CDP-24-004-C2EMPI project

\medskip
Received XXXXXX; revised XXXXXX.
\medskip

\end{document}